\documentclass{article}
\usepackage{latexsym}
\usepackage{amstext}
\usepackage{amsmath}
\usepackage{amssymb}
\usepackage{amsthm}

\newtheorem{theorem}{Theorem}
\newtheorem{lemma}{Lemma}
\newtheorem{corollary}{Corollary}
\theoremstyle{definition}

\newtheorem{example}{Example}

\def \co {\mathcal{O}}
\def \kbar {\bar{k}}
\def \bd {d}
\def \cD {\mathcal{D}}
\def \cC {\mathcal{C}}
\def \Q {\mathbb{Q}}

\def \P {\mathbb{P}}
\def \A {\mathbb{A}}
\def \C {\mathbb{C}}

\DeclareMathOperator{\Tr}{Tr}

\DeclareMathOperator{\Spec}{Spec}
\DeclareMathOperator{\Supp}{Supp}

\begin{document}
\bibliographystyle{amsplain}
\title{Vojta's Inequality and Rational and Integral Points of Bounded Degree on Curves}
\author{Aaron Levin\\adlevin@math.brown.edu}
\maketitle
\begin{abstract}
Let $C\subset C_1\times C_2$ be a curve of type $(d_1,d_2)$ in the product of the two curves $C_1$ and $C_2$.  Let $\nu$ be a positive integer.  We prove that if a certain inequality involving $d_1$, $d_2$, $\nu$, and the genera of the curves $C_1$, $C_2$, and $C$ is satisfied, then the set of points $\{P\in C(\kbar)\mid [k(P):k]\leq \nu \}$ is finite for any number field $k$.  We prove a similar result for integral points of bounded degree on $C$.  These results are obtained as consequences of an inequality of Vojta which generalizes the Roth-Wirsing theorem to curves.
\end{abstract}
\section{Introduction}
In \cite{Vo1}, Vojta proved the following theorem.
\begin{theorem}[(Vojta)]
\label{V}
Let $C$ be a nonsingular curve defined over a number field $k$.  Let $X$ be a regular model for $C$ over the ring of integers of $k$.  Let $K$ be the canonical divisor of $C$, $A$ an ample divisor on $C$, and $D$ an effective divisor on $C$ without multiple components.  Let $S$ be a finite set of places of $k$.  Let $\nu$ be a positive integer and let $\epsilon>0$.  Then
\begin{equation}
\label{Vojta}
m_S(D,P)+h_K(P)\leq d_a(P)+\epsilon h_A(P)+O(1)
\end{equation}
for all points $P\in C(\kbar)\setminus \Supp D$ with $[k(P):k]\leq \nu$.
\end{theorem}
Here $h_D$ is a logarithmic height associated to the divisor $D$, $m_S(D,P)$ is a proximity function, and $d_a(P)$ is the arithmetic discriminant of \cite{Vo2}, whose definition we recall below.  We refer the reader to \cite{La}, \cite{Vo3}, and \cite{Vo1} for definitions and properties of heights and proximity functions.

The inequality (\ref{Vojta}) is a vast generalization of the theorems of Roth and Wirsing. In particular, it implies Faltings' theorem (Mordell's conjecture).  As a consequence of (\ref{Vojta}), Song and Tucker \cite{Tu} show
\begin{corollary}[(Song, Tucker, Vojta)]
\label{co1}
Let $C$ and $C'$ be nonsingular curves of genus $g$ and $g'$, respectively, defined over a number field $k$.  Let $\phi:C\to C'$ be a dominant $k$-morphism.  If
\begin{equation}
\label{ineq}
g-1>(\nu+g'-1)\deg \phi
\end{equation}
for some positive integer $\nu$, then the set
\begin{equation}
\label{set}
\{P\in C(\kbar)\mid [k(P):k]\leq \nu \text{ and } k(\phi(P))=k(P)\}
\end{equation}
is finite.
\end{corollary}
Vojta noted the case $C'=\mathbb{P}^1$ of the corollary.  Note that the condition $k(\phi(P))=k(P)$ in Theorem \ref{co1} precludes one from deducing a finiteness result on algebraic points with $[k(P):k]\leq \nu$.  Of course, this condition in the theorem is necessary (consider, for example, hyperelliptic curves of genus $g>3$).  If we are given more than one dominant morphism of $C$ to a curve where (\ref{ineq}) holds, it is natural to try to prove a finiteness result without the $k(\phi(P))=k(P)$ condition in (\ref{set}).  Clearly we need the two maps to be independent in some sense.  More precisely, we will assume that we are given a birational morphism of $C$ into a product of curves.  In addition to rational points, we will study integral points on $C$.

Let $S$ be a finite set of places of $k$ and let $\co_{k,S}$ be the ring of $S$-integers of $k$.  Let $D$ be an effective divisor on $C$.  If $D\neq 0$, we call a set $T\subset C(\kbar)\setminus \Supp D$ a set of $(D,S)$-integral points on $C$ if there exists an affine embedding $C\setminus \Supp D\subset \A^m$ such that every point $P\in T$ has $S$-integral coordinates, i.e., each coordinate of $P$ in $\A^m$ lies in the integral closure of $\co_{k,S}$ in $\kbar$.  If $D=0$, then we call any subset of $C(\kbar)$ a set of $(D,S)$-integral points.  Our main theorem is
\begin{theorem}
\label{prod}
Let $C$, $C_1$, and $C_2$ be nonsingular curves of genus $g$, $g_1$, and $g_2$, respectively, all defined over a number field $k$.  Let $S$ be a finite set of places $k$.  Let $\phi:C\to C_1\times C_2$ be a birational morphism.  Let $\pi_1$ and $\pi_2$ be the projections of $C_1\times C_2$ onto the first and second factors, respectively.  Suppose that $\pi_1\circ \phi$ and $\pi_2\circ \phi$ are dominant morphisms and let $d_1=\deg \pi_1\circ \phi$ and $d_2=\deg \pi_2\circ \phi$.  Let $D=\sum_{i=1}^r P_i$ be an effective divisor on $C$, defined over $k$, with $P_1,\ldots, P_r$ distinct points of $C(\kbar)$.  If
\begin{equation}
\label{ineq2}
2g-2+r>\max\{(\nu+g_1-1)2d_1,(\nu+g_2-1)2d_2,(\nu+2g_1-2)d_1+(\nu+2g_2-2)d_2\}
\end{equation}
for some positive integer $\nu$, then any set of $(D,S)$-integral points
\begin{equation*}
T\subset \{P\in C(\kbar)\mid [k(P):k]\leq \nu \}
\end{equation*}
is finite.  In particular, if \eqref{ineq2} is satisfied with $r=0$, then the set
\begin{equation*}
\{P\in C(\kbar)\mid [k(P):k]\leq \nu \}
\end{equation*}
is finite.
\end{theorem}
\section{Some Examples and Corollaries}
We first give two examples which show that the inequality (\ref{ineq2}) is sharp in the sense that Theorem \ref{prod} is false if ``$>$" is replaced by ``$\geq$" in (\ref{ineq2}).
\begin{example}
Let $C$ be a nonsingular curve, defined over a number field $k$, of bidegree $(d_1,d_2)$ on $C_1\times C_2=\P^1\times \P^1$ with $d_1\geq d_2 >0$.  Let $P,Q\in \P^1(k)$ be two points above which $\phi_2$ is unramified, and let $D=P+Q$.  Over sufficiently large number fields $k$, there are infinitely many $k$-rational $(D,S)$-integral points on $\P^1$.  Pulling back these integral points by $\phi_2$, we obtain infinitely many $(\phi_2^*D,S)$-integral points on $\P^1\times \P^1$ (of degree $\leq d_2=\nu$ over $k$), where $\phi_2^*D$ is a sum of $r=2d_2$ distinct points.  We have $g=(d_1-1)(d_2-1)$ and we see that equality holds in (\ref{ineq2}).  
\end{example}
\begin{example}
Let $C_1\times C_2=\P^1\times E$, where $E$ is an elliptic curve defined over a number field $k$.  Let $d_1>d_2+1>2$.  Let $C$ be a nonsingular curve, defined over a number field $k$, of type $(d_1,d_2)$ on $\P^1\times E$ (i.e., $\deg \pi_1|_C=d_1$ and $\deg \pi_2|_C=d_2$).  Then by the adjunction formula, $g=g(C)=d_1(d_2-1)+1$.  Let $\nu=d_2$ and $r=0$.  Then a simple calculation shows that equality is achieved in (\ref{ineq2}), but the set $\{P\in C(\kbar)\mid [k(P):k]\leq \nu\}$ is infinite for sufficiently large $k$ since $C$ has a degree $\nu=d_2$ map down to $E$.
\end{example}
Note that when $C_1\times C_2=\P^1\times \P^1$, the inequality (\ref{ineq2}) simplifies to 
\begin{equation*}
2g-2+r>\max\{2(\nu-1)d_1,2(\nu-1)d_2\}.
\end{equation*}
As a curve of degree $d$ in $\mathbb{P}^2$ can be mapped birationally onto a curve of bidegree $(d-1,d-1)$ in $\P^1\times \P^1$, we obtain
\begin{corollary}
\label{co2}
Let $C\subset \P^2$ be a curve, defined over a number field $k$, of degree $d$ and geometric genus $g$.  Let $S$ be a finite set of places of $k$.  Let $D=\sum_{i=1}^r P_i$ be an effective divisor on $C$, defined over $k$, with $P_1,\ldots, P_r$ distinct points of $C(\kbar)$.  If
\begin{equation}
\label{ineq3}
2g-2+r>2(\nu-1)(d-1)
\end{equation}
for some positive integer $\nu$, then any set of $(D,S)$-integral points
\begin{equation*}
T\subset \{P\in C(\kbar)\mid [k(P):k]\leq \nu \}
\end{equation*}
is finite.  In particular, if $g-1>(\nu-1)(d-1)$ then the set
\begin{equation*}
\{P\in C(\kbar)\mid [k(P):k]\leq \nu \}
\end{equation*}
is finite.
\end{corollary}
By definition, the geometric genus of $C$ is the genus of the normalization of $C$.  For nonsingular plane curves, a better theorem on rational points has been proven by Debarre and Klassen \cite{De2} using Falting's theorem on rational points on subvarieties of abelian varieties.
\begin{theorem}[(Debarre, Klassen)]
Let $C\subset \P^2$ be a nonsingular curve of degree $d$, defined over a number field $k$, that does not admit a map of degree $\leq d-2$ onto a genus one curve (this is automatically satisfied if $d\geq 7$).  Then the set 
\begin{equation*}
\{P\in C(\kbar)\mid [k(P):k]\leq d-2 \}
\end{equation*}
is finite.
\end{theorem}
Recall that a curve is called hyperelliptic (respectively bielliptic) if it admits a map of degree two onto a curve of geometric genus zero (respectively one).
Harris and Silverman \cite{HS} have shown (again using Falting's theorem on subvarieties of abelian varieties)
\begin{theorem}[(Harris, Silverman)]
\label{HS}
Let $C$ be a nonsingular curve defined over a number field $k$.  If $C$ is not hyperelliptic or bielliptic then the set $\{P\in C(\kbar)\mid [k(P):k]\leq 2 \}$ is finite.
\end{theorem}
A similar theorem is true for degree three rational points (see \cite{AH}), but not for degrees four and higher (see \cite{De}).  Similarly, for integral points, Corvaja and Zannier \cite{CZ} have shown 
\begin{theorem}[(Corvaja, Zannier)]
\label{CZ}
Let $C$ be a nonsingular curve defined over a number field $k$.  Let $S$ be a finite set of places of $k$.  Let $D=\sum_{i=1}^r P_i$ be an effective divisor on $C$, defined over $k$, with $P_1,\ldots, P_r$ distinct points of $C(\kbar)$.  Let $T\subset \{P\in C(\kbar)\mid [k(P):k]\leq 2 \}$ be a set of $(D,S)$-integral points.  Then
\begin{enumerate}
\item  If $r>4$ then $T$ is finite.
\item  If $r>3$ and $C$ is not hyperelliptic then $T$ is finite.
\end{enumerate}
\end{theorem}
Additionally, in the case $C$ is hyperelliptic and $r=4$ (where $T$ may be infinite), Corvaja and Zannier show how to parametrize all but finitely many of the quadratic integral points.
The proof of Theorem \ref{CZ} in \cite{CZ} makes use of an appropriate version of the Schmidt subspace theorem.
We now show that Corollary \ref{co2} implies a slight improvement to this theorem.  Specifically, we show that the inequality in part (b) can be improved to cover the case $r=3$.
\begin{theorem}
\label{CZ2}
Let $C$ be a nonsingular curve defined over a number field $k$.  Let $S$ be a finite set of places of $k$.  Let $D=\sum_{i=1}^r P_i$ be an effective divisor on $C$, defined over $k$, with $P_1,\ldots, P_r$ distinct points of $C(\kbar)$.  Let $T\subset \{P\in C(\kbar)\mid [k(P):k]\leq 2 \}$ be a set of $(D,S)$-integral points.  Then
\begin{enumerate}
\item  If $r>4$ then $T$ is finite.
\item  If $r>2$ and $C$ is not hyperelliptic then $T$ is finite.
\end{enumerate}
\end{theorem}
\begin{proof}
By Corollary \ref{co2}, to prove (a) it suffices to show that any curve $C$ of genus $g$ has a birational plane model of degree $g+2$.  Since any divisor of degree $2g+1$ on $C$ is very ample and nonspecial, we obtain an embedding of $C$ as a degree $2g+1$ curve in $\P^{g+1}$.  Projecting from the linear span of $g-1$ general points of $C$, we obtain a birational map $\phi:C\to \P^2$ with $\deg \phi(C)=g+2$ (see \cite[p. 109]{ACGH}).

Similarly, to prove (b) it suffices to show that if $C$ has genus $g$ and is not hyperelliptic, then $C$ has a birational  plane model of degree $g+1$.  Since $C$ is not hyperelliptic, the canonical embedding realizes $C$ as a curve of degree $2g-2$ in $\P^{g-1}$.  Projecting from the linear span of $g-3$ general points of $C$, we obtain a plane curve of degree $g+1$ birational to $C$.
\end{proof}
As noted in \cite{CZ}, Vojta's conjecture predicts that the inequality in (b) can be improved to $r>0$.  It is unclear to what extent this follows from Theorem \ref{prod}.  For instance, Theorem \ref{prod} implies that one may take $r>0$ in Theorem \ref{CZ2} for any nonsingular bielliptic curve $C$ of type $(a,2)$, $a>3$, on $\P^1\times E$ (of course, by Theorem \ref{HS}, we need only consider bielliptic curves in (b)).  
\section{Proofs of Results}
Let $C$ be a nonsingular curve defined over a number field $k$.  Let $R$ denote the ring of integers of $k$ and let $B=\Spec R$.  Let $\pi:X\to B$ be a regular model for $C$ over $R$.  For every complex embedding $\sigma:k\hookrightarrow \C$ we have a canonical volume form on $C_\sigma=C\times_\sigma \C$ and an associated canonical Green's function $g_\sigma$.  With this data one can define intersections of Arakelov divisors (see \cite{La2}).  Let $P\in C(\kbar)$ and let $E_P$ denote the horizontal prime divisor on $X$ corresponding to $P$ (we will also denote the curve on $X$ corresponding to $P$ by $E_P$).  Let $\omega_{X/B}$ denote the relative dualizing sheaf, with its canonical Arakelov metric \cite[Ch. 4]{La2}.  We then define the arithmetic discriminant $d_a(P)$ by
\begin{equation*}
d_a(P)=\frac{(E_P.(\omega_{X/B}+E_P))}{[k(P):\Q]}.
\end{equation*}
Of course, contrary to the notation, $d_a(P)$ depends on more data than just $P$.  We can also give an alternative formula for $d_a(P)$.  Let $L=k(P)$.  Then $E_P=\Spec A$, where $A$ is an order of the number field $L$.  Let
\begin{equation*}
W_{A/R}=\{b\in L\mid \Tr_{L/k}(bA)\subset R\}
\end{equation*}
be the Dedekind complementary module.  It is a fractional ideal of $A$ containing $A$.  For a fractional ideal $\mathfrak{a}$ of  $A$, we define the fractional ideal
\begin{equation*}
\mathfrak{a}^{-1}=\{x\in L\mid x\mathfrak{a}\subset  A\}.
\end{equation*}
In arbitrary orders, one may not necessarily have $\mathfrak{a}\mathfrak{a}^{-1}=A$.  We now define the Dedekind different (of $A$ over $R$) as
\begin{equation*}
\cD_{A/R}=W_{A/R}^{-1}.
\end{equation*}
This is an integral ideal of $A$.  For a nice discussion of the relation between the different, discriminant, and conductor of an order, we refer the reader to the article by Del Corso and Dvornicich \cite{Dv}.  Now define
\begin{equation*}
\bd_{A/R}=\frac{\log [A:\cD_{A/R}]}{[L:\Q]}
\end{equation*}
Let $S_\infty$ be the set of archimedean places of $k$ and let $v\in S_\infty$.  Let 
\begin{equation*}
E_v=E_P\times \C_v=\{P_{v,1},\ldots,P_{v,[L:k]}\}
\end{equation*}
be the points in $C_v=C\times \C_v$ into which $E_P$ splits.
By the Arakelov adjunction formula \cite[Th. 5.3]{La2}, we have
\begin{equation}
\label{ar}
d_a(P)=d_{A/R}+\frac{1}{[L:\Q]}\sum_{v\in S_\infty}\sum_{i \neq j}N_v\lambda_v(P_{v,i},P_{v,j})
\end{equation}
where $N_v=[k_v:\Q_v]$, and $\lambda_v=\frac{1}{2}g_v$ (with $g_v$ normalized as in \cite{La2}).  We will use that $\lambda_v$ is a Weil function for the diagonal $\Delta_v$ in $C_v\times C_v$, i.e., if the Cartier divisor $\Delta_v$ is locally represented by a function $f$ on the open set $U$, then there exists a continuous function $\alpha$ on $U$ such that 
\begin{equation*}
\lambda_v(P)=-\log |f(P)|+\alpha(P)
\end{equation*}
for all $P\in U\setminus \Delta_v$.
\begin{theorem}
\label{main}
Let $C_1$, $C_2$, and $C_3$ be nonsingular curves defined over $k$ and let $X_1$, $X_2$, and $X_3$ be regular models over $R$ for the respective curves.  Let $\phi:X_3\to X_1 \times X_2$ be a birational morphism onto its image.  Let $\phi_1$ and $\phi_2$ denote $\phi$ composed with the projection map of $X_1\times X_2$ onto the first and second factor, respectively.  Let $P\in C_3(\kbar)$.  Then
\begin{equation}
\label{triangle2}
d_a(P)\leq d_a(\phi_1(P))+d_a(\phi_2(P))+O(1).
\end{equation}
\end{theorem}
Our strategy is to break up $d_a$ into a finite and infinite part as in (\ref{ar}), and then prove the inequality for each part separately.  Since there is an $O(1)$ term, we can clearly ignore the finite set $Z$ of $C(\kbar)$ on which $\phi_{\kbar}$ fails to be invertible.  To prove the inequality for the finite part, $d_{A/R}$, of (\ref{ar}), we use the following lemma.
\begin{lemma}
\label{l1}
Let $R$ be the ring of integers of a number field $k$.  Let $A_1$ and $A_2$ be $R$-orders of the number fields $L_1$ and $L_2$, respectively (with some fixed embedding in $\kbar$).  Let $L_3=L_1L_2$ and let $A_3=A_1A_2$.  If $A_1$, $A_2$, and $A_3$ are Gorenstein rings then 
\begin{equation}
\label{triangle}
\bd_{A_3/R}\leq \bd_{A_1/R}+\bd_{A_2/R}
\end{equation}
\end{lemma}
\begin{proof}
As shown in \cite{Dv}, an $R$-order $A$ is Gorenstein if and only if $\cD_{A/R}$ is an invertible ideal of $A$ (see \cite{Ba} for the many equivalent definitions of a Gorenstein ring).  Let $A_i'$ denote the integral closure of $A_i$ in $L_i$ for $i=1,2,3$.  For the Gorenstein rings $A_1$, $A_2$, and $A_3$ we have the relations (see \cite[Prop. 3]{Dv})
\begin{equation}
\label{cond}
\cD_{A_i/R}A_i'=\cC_{A_i} \cD_{A_i'/R}, \quad i=1,2,3,
\end{equation}
where
\begin{equation*}
\cC_{A_i}=\{x\in A_i'\mid xA_i'\subset A_i\}
\end{equation*}
is the conductor of $A_i$.
For an invertible ideal $\mathfrak{a}$ of ${A_3}$ (see \cite[Th. 1]{Dv}), 
\begin{equation*}
[{A_3}:\mathfrak{a}]=[{A_3}':\mathfrak{a}{A_3}'].
\end{equation*}
Now to prove the lemma, it suffices to show that 
\begin{equation*}
\cD_{A_1/R}\cD_{{A_2}/R}{A_3}'\subset \cD_{{A_3}/R}{A_3}'.
\end{equation*}
Indeed, this inclusion gives 
\begin{equation*}
[{A_3}:\cD_{{A_3}/R}]=[{A_3}':\cD_{{A_3}/R}{A_3}']\leq [{A_3}':\cD_{A_1/R}\cD_{{A_2}/R}{A_3}']
\end{equation*}
which is equivalent to (\ref{triangle}) as
\begin{align*}
[{A_3}':\cD_{A_1/R}\cD_{{A_2}/R}{A_3}']&=[{A_3}':\cD_{A_1/R}{A_3}'][{A_3}':\cD_{{A_2}/R}{A_3}']\\
&=[A_1':\cD_{A_1/R}A_1']^{[L_3:L_1]}[{A_2}':\cD_{{A_2}/R}{A_2}']^{[L_3:L_2]}\\
&=[A_1:\cD_{A_1/R}]^{[L_3:L_1]}[{A_2}:\cD_{{A_2}/R}]^{[L_3:L_2]}.
\end{align*}
We now show that $\cD_{A_1/R}\cD_{{A_2}/R}{A_3}'\subset \cD_{{A_3}/R}{A_3}'$.  By (\ref{cond}),
\begin{equation*}
\cD_{A_1/R}\cD_{{A_2}/R}{A_3}'=\cC_{A_1} \cD_{A_1'/R}\cC_{A_2} \cD_{{A_2}'/R}{A_3}'
\end{equation*}
and
\begin{equation*}
\cD_{{A_3}/R}{A_3}'=\cC_{A_3} \cD_{{A_3}'/R}=\cC_{A_3} \cD_{{A_3}'/A_1'}\cD_{A_1'/R}.
\end{equation*}
Therefore we need to show $\cC_{A_1}\cC_{A_2} \cD_{{A_2}'/R}{A_3}'\subset \cC_{A_3} \cD_{{A_3}'/A_1'}$.
It is a standard fact that $\cD_{{A_2}'/R}$ is generated by elements of the form $f'(\alpha)$, where $\alpha \in {A_2}'$, $k(\alpha)=L_2$, and $f$ is the minimal polynomial of $\alpha$ over $k$.  Let $g$ be the minimal polynomial of $\alpha$ over $L_1$.  Note that $L_1(\alpha)=L_3$ and that $g'(\alpha)$ divides $f'(\alpha)$ in ${A_3}'$.  It is easily shown that $g'(\alpha){A_3}'=\cC_{A_1'[\alpha]}\cD_{{A_3}'/A_1'}$.  We have 
\begin{equation*}
\cC_{A_1}\cC_{A_2}\cC_{A_1'[\alpha]}\subset \cC_{A_3}
\end{equation*}
since 
\begin{equation*}
\cC_{A_1}\cC_{A_2}\cC_{A_1'[\alpha]}{A_3}'\subset\cC_{A_1}\cC_{A_2} A_1'[\alpha]\subset \cC_{A_1}\cC_{A_2}A_1'{A_2}'\subset A_1{A_2}={A_3}.
\end{equation*}
Therefore
\begin{equation*} 
\cC_{A_1}\cC_{A_2}f'(\alpha)\subset \cC_{A_1}\cC_{A_2}\cC_{A_1'[\alpha]}\cD_{{A_3}'/A_1'}\subset \cC_{A_3} \cD_{{A_3}'/A_1'}
\end{equation*}
As $D_{{A_2}'/R}$ was generated by the $f'(\alpha)$, we obtain $\cC_{A_1}\cC_{A_2} \cD_{{A_2}'/R}{A_3}'\subset \cC_{A_3} \cD_{{A_3}'/A_1'}$ as desired.
\end{proof}
Now let $E_P=E_3=\Spec A_3$ be the prime horizontal divisor corresponding to $P\in C(\kbar)\setminus Z$, and let $\phi_1(E_P)=E_1=\Spec A_1$ and $\phi_2(E_P)=E_2=\Spec A_2$.  Note that $A_1$ and $A_2$ are naturally subrings of $A_3$ (via $\phi_1$ and $\phi_2$) and $A_3=A_1A_2$.  Indeed, the closed immersion $\phi:E_P\to X_1\times X_2$ factors through $E_1\times E_2$, and therefore the natural map $A_1 \otimes A_2\to A_3$ is surjective.  Since $X_1$, $X_2$, and $X_3$ were assumed regular, $E_P$, $E_1$, and $E_2$ are locally complete intersections (they are Cartier divisors).  This implies in particular that $A_1$, $A_2$, and $A_3$ are Gorenstein rings.  Therefore, using Lemma \ref{l1}, we have proved the finite part of the inequality (\ref{triangle2}), i.e., the inequality (\ref{triangle}).  

We now consider the archimedean part of (\ref{triangle2}).  With notation as above, let $L_1$, $L_2$, and $L_3$ be the quotient fields of $A_1$, $A_2$, and $A_3$.  Let $v \in S_\infty$.  Let $E_{iv}$ be the set of points of $E_i\times \C_v$,  $i=1,2,3$.  Let $\lambda_{\Delta_1}$, $\lambda_{\Delta_2}$, and $\lambda_{\Delta_3}$ denote the Weil functions (relative to $v$) of (\ref{ar}) for $C_1$, $C_2$, and $C_3$, respectively.  Here $\Delta_i$ is the diagonal of $C_{iv}\times C_{iv}$.  Then it suffices to prove
\begin{lemma}
\label{l2}
In the notation above,
\begin{multline*}
\frac{1}{[L_3:\Q]}\sum_{\substack{P,Q\in E_{3v}\\ P\neq Q}}\lambda_{\Delta_3}(P,Q)\leq \frac{1}{[L_1:\Q]}\sum_{\substack{P,Q\in E_{1v}\\ P\neq Q}}\lambda_{\Delta_1}(P,Q) +\\
\frac{1}{[L_2:\Q]}\sum_{\substack{P,Q\in E_{2v}\\ P\neq Q}}\lambda_{\Delta_2}(P,Q)+O(1).
\end{multline*}
\end{lemma}
The lemma will follow easily from the following ``distribution relation" of Silverman \cite[Prop. 6.2(b)]{Si} (proved by Silverman in greater generality).
\begin{theorem}[(Silverman)]
\label{Sil}
Let $C$ and $C'$ be nonsingular complex curves.  Let $\phi:C\to C'$ be a morphism.  Let $\Delta$ and $\Delta'$ denote the diagonals of $C\times C$ and $C'\times C'$, respectively.  Let $\lambda_{\Delta}$ and $\lambda_{\Delta'}$ be Weil functions associated to $\Delta$ and $\Delta'$ (under the usual complex absolute value).  Then for any $P\in C$ and $q\in C'$ with $\phi(P)\neq q$, 
\begin{equation*}
\lambda_{\Delta'}(\phi(P),q)=\sum_{Q\in \phi^{-1}(q)}e_\phi(Q/q)\lambda_{\Delta}(P,Q)+O(1) 
\end{equation*}
where $e_\phi(Q/q)$ is the ramification index of $\phi$ at $Q$.
\end{theorem}
\begin{proof}[Proof of Lemma \ref{l2}]
Denote by $\phi$, $\phi_1$, and $\phi_2$ the same maps base extended to $C_{3v}$.  Let $P,Q \in E_{3v}$, $P\neq Q$.  Since we assumed $P\notin Z$, either $\phi_1(P)\neq \phi_1(Q)$ or $\phi_2(P)\neq \phi_2(Q)$.  Note also that the maps $E_{3v}\to E_{1v}$ and $E_{3v}\to E_{2v}$ are $[L_3:L_1]$-to-$1$ and $[L_3:L_2]$-to-$1$ maps respectively.  Thus we obtain (modulo bounded functions independent of $E_{1v}$, $E_{2v}$, and $E_{3v}$)
\begin{align*}
\sum_{\substack{P,Q\in E_{3v}\\ P\neq Q}}\lambda_{\Delta_3}(P,Q)&\leq \sum_{\substack{P,Q\in E_{3v}\\ \phi_1(P)\neq \phi_1(Q)}}\lambda_{\Delta_3}(P,Q)+\sum_{\substack{P,Q\in E_{3v}\\ \phi_2(P)\neq \phi_2(Q)}}\lambda_{\Delta_3}(P,Q)\\
&\leq \sum_{P\in E_{3v}}\sum_{\substack{Q\in \phi_1^{-1}(E_{1v})\\ \phi_1(P)\neq \phi_1(Q)}}\lambda_{\Delta_3}(P,Q)+ \sum_{P\in E_{3v}}\sum_{\substack{Q\in \phi_2^{-1}(E_{2v})\\ \phi_2(P)\neq \phi_2(Q)}}\lambda_{\Delta_3}(P,Q)\\
&\leq \sum_{P\in E_{3v}}\sum_{\substack{q\in E_{1v}\\ \phi_1(P)\neq q}}\lambda_{\Delta_1}(\phi_1(P),q)+\sum_{P\in E_{3v}}\sum_{\substack{q\in E_{2v}\\ \phi_1(P)\neq q}}\lambda_{\Delta_2}(\phi_2(P),q)\\
&\leq [L_3:L_1]\sum_{\substack{p,q\in E_{1v}\\ p\neq q}}\lambda_{\Delta_1}(p,q)+[L_3:L_2]\sum_{\substack{p,q\in E_{2v}\\ p\neq q}}\lambda_{\Delta_2}(p,q)
\end{align*}
Dividing everything by $[L_3:\Q]$ gives the lemma.
\end{proof}
Theorem \ref{main} now follows from Lemma \ref{l1} and Lemma \ref{l2}.  We now prove Theorem~\ref{prod} from the introduction.  We will need the following estimate of Song and Tucker (see \cite{Tu2} and \cite{Tu}) for $d_a(P)$ on a curve.
\begin{lemma}[(Song, Tucker)]
\label{Tu}
Let $C$ be a nonsingular curve defined over a number field $k$ with canonical divisor $K$.  Let $X$ be a regular model for $C$ over the ring of integers of $k$.  Let $A$ be an ample divisor on $C$ and let $\epsilon>0$.  Then
\begin{equation*}
d_a(P)\leq h_K(P)+(2[k(P):k]+\epsilon)h_A(P)+O([k(P):k]).
\end{equation*}
\end{lemma}
\begin{proof}[Proof of Theorem \ref{prod}]
Let $T$ be as in the hypotheses of Theorem \ref{prod} and suppose that the inequality (\ref{ineq2}) of Theorem \ref{prod} is satisfied.
Consider the three sets 
\begin{align*}
&T_1=\{P\in T\mid [k(\phi_1(P)):k]=[k(P):k]\}\\
&T_2=\{P\in T\mid [k(\phi_2(P)):k]=[k(P):k]\}\\
&T_3=\{P\in T\mid [k(\phi_1(P)):k]<[k(P):k],[k(\phi_2(P)):k]<[k(P):k]\}.
\end{align*}
Clearly $T=T_1\cup T_2\cup T_3$.  As we assumed $2g-2+r>(\nu+g_1-1)2d_1$ and $2g-2+r>(\nu+g_2-1)2d_2$, it follows from a trivial generalization of Corollary \ref{co1} that $T_1$ and $T_2$ are finite.  So we are reduced to showing that if $2g-2+r>(\nu+2g_1-2)d_1+(\nu+2g_2-2)d_2$ then $T_3$ is finite.  Let $K$, $K_1$, and $K_2$ denote the canonical divisors of $C$, $C_1$, and $C_2$, respectively.  Let $h$, $h_1$, and $h_2$ denote heights associated to some degree one divisor on $C$, $C_1$ and $C_2$, respectively.  Using Theorem \ref{V}, Theorem \ref{main}, and Lemma \ref{Tu}, we get, for any $\epsilon>0$,
\begin{align*}
m_S(D,P)+h_K(P)&\leq d_a(P)+\epsilon h(P)+O(1)\\
&\leq d_a(\phi_1(P))+d_a(\phi_2(P))+\epsilon h(P)+O(1)\\
&\leq h_{K_1}(\phi_1(P))+(2[k(\phi_1(P)):k]+\epsilon)h_1(\phi_1(P))+\\
& \qquad \qquad h_{K_2}(\phi_2(P))+(2[k(\phi_2(P)):k]+\epsilon)h_2(\phi_2(P))+O(1).
\end{align*}
Note that for $P\in T_3$, $[k(\phi_1(P)):k]\leq \nu/2$ and $[k(\phi_2(P)):k]\leq \nu/2$, since $k(\phi_1(P))$ and $k(\phi_2(P))$ are both proper subfields of $k(P)$.  Since $T$ is a set of $(D,S)$-integral points, $m_S(D,P)=h_D(P)+O(1)$ for $P\in T$.  Using functoriality of heights and quasi-equivalence of heights associated to numerically equivalent divisors, we obtain, for any $\epsilon>0$,
\begin{align*}
(2g-2+r)h(P)\leq ((\nu+2g_1-2)d_1+(\nu+2g_2-2)d_2+\epsilon)h(P)+O(1)
\end{align*}
for $P\in T_3$.  Taking $\epsilon<1$, since there are only finitely many points of bounded degree and bounded height, we see that if
\begin{equation*}
2g-2+r>(\nu+2g_1-2)d_1+(\nu+2g_2-2)d_2
\end{equation*}
then $T_3$, and hence $T$, must be finite.
\end{proof}
{\bf Acknowledgments.}  I would like to thank Joe Silverman for helpful conversations and for pointing me to the reference for Theorem \ref{Sil}.
\bibliography{Arakelov}
\end{document}